\documentclass{amsart}
\usepackage{amsfonts}
\usepackage{amsmath,amssymb}
\usepackage{amsthm}
\usepackage{amscd}
\usepackage{graphics}
\usepackage{graphicx}

\theoremstyle{remark}{
\newtheorem{Def}{{\rm Definition}}

\newtheorem{Rem}{{\rm Remark}}

}
\theoremstyle{plain}
{
\newtheorem{Cor}{Corollary}

\newtheorem{MainThm}{Main Theorem}

}

\begin{document}
\title[Regions surrounded by circles whose Poincar\'e-Reeb graphs are trees]{Regions surrounded by circles whose Poincar\'e-Reeb graphs are trees}
\author{Naoki kitazawa}

\address{Institute of Mathematics for Industry, Kyushu University, 744 Motooka, Nishi-ku Fukuoka 819-0395, Japan\\
 TEL (Office): +81-92-802-4402 \\
 FAX (Office): +81-92-802-4405 \\
}
\email{naokikitazawa.formath@gmail.com}
\urladdr{https://naokikitazawa.github.io/NaokiKitazawa.html}
\maketitle
\begin{abstract}
Regions in the Euclidean plane surrounded by circles are fundamental geometric and combinatorial objects. Related studies have been done and we cannot explain them precisely, or roughly, well.

We study such regions whose {\it Poincar\'e-Reeb graphs} are trees and investigate the trees obtained by a certain inductive rule from a disk in the plane. The {\it Poincar\'e-Reeb} graph of such a region is a graph whose underlying set is the set of all components of level sets of the restriction of the canonical projection to the closure and whose vertices are points corresponding to the components containing {\it singular} points. Related studies were started by the author, motivated by importance and difficulty of explicit construction of a real algebraic map onto a prescribed closed region in the plane.

\end{abstract}
\section{Introduction.}
\label{sec:1}
Regions in the Euclidean plane surrounded by circles are fundamental geometric and combinatorial objects. Related studies have been done, where we cannot explain them precisely, or roughly, well. Here, refer to the preprints of the author \cite{kitazawa4, kitazawa5, kitazawa6} with \cite{kitazawa3}, mainly. We do not assume related knowledge. We study the shape of the region. We can define its {\it Poincar\'e-Reeb graph}. The underlying set consists of all components of level sets of the restriction of the canonical projection to the closure and a point is a vertex if and only if it corresponds to a component containing {\t singular} points of the region. We define related notions rigorously, later.

The author has been interested in reconstructing a real algebraic map locally like a moment-map with prescribed image being also a region in the plane. See the paper \cite{kitazawa2} and \cite{kitazawa3} mainly. This is also related to the difficulty of construction of real algebraic objects, which is different from the existence and approximation theory by \cite{nash, tognoli}, surveyed in \cite{kollar} for example. Such objects exist plentifully.

Related to this, the author has formulated a rule of obtaining a new region by adding a small circle centered at a point in the boundary of the given region and another rule, adding a circle passing a small chord connecting two points in the boundary of the region and sufficiently close to the chord, in \cite{kitazawa3}. The preprint concentrates on reconstruction of real algebraic maps and the rules are formulated as tools for inductive construction. \cite{kitazawa3} are studies on these rules and types of changes of Poincar\'e-Reeb graphs are explicitly studied. There cases which are not discussed in the original study \cite{kitazawa3} are also studied.
\subsection{Fundamental terminologies and notation.}
Let $\emptyset$ denote the empty set. For a topological space $X$ and its subspace $Y$, let ${\overline{Y}}^X$ denote the closure. Let $X$ be a topological space having the structure of some cell complex the dimensions of cell of which are bounded. We can define the dimension $\dim X$ uniquely as a non-negative integer and this gives a topological invariant. 
A topological manifold is well-known to be homeomorphic to a CW complex. A smooth manifold is regarded as a polyhedron and has the structure of a certain polyhedron canonically ({\rm PL manifold}). A topological space with the structure of a polyhedron whose dimension is at most $2$ has the unique structure of a polyhedron. This holds for topological manifolds of dimension at most $3$. This is presented in \cite{moise}. We use ${\rm Int}\ X$ for the interior of a manifold $X$ and $\partial X:=X-{\rm Int}\ X$ denotes the boundary.
 
Let ${\mathbb{R}}^k$ denote the $k$-dimensional Euclidean space, a simplest smooth manifold, and the Riemannian manifold equipped with the so-called standard Euclidean metric. We use $\mathbb{R}:={\mathbb{R}}^1$.
For a point $x \in {\mathbb{R}}^k$, let $||x|| \geq 0$ denote the distance between $x$ and the origin $0 \in {\mathbb{R}}^k$.  Let $S^k:=\{x \in {\mathbb{R}}^{k+1} \mid ||x||=1\}$ denote the $k$-dimensional {\it unit sphere}, a $k$-dimensional smooth compact submanifold of ${\mathbb{R}}^{k+1}$ with no boundary. It is connected for $k \geq 1$ and a discrete two-point set for $k=0$. It is the zero set of the polynomial $||x||^2-1={\Sigma}_{j=1}^{k+1} {x_j}^2-1$ with $x:=(x_1,\cdots,x_{k+1})$. Let $D^k:=\{x \in {\mathbb{R}}^{k} \mid ||x|| \leq 1\}$ denote the $k$-dimensional {\it unit disk}, a $k$-dimensional smooth compact and connected submanifold of ${\mathbb{R}}^{k}$. We also have $\partial D^k=S^{k-1}$.

For a differentiable manifold $X$, let $T_x X$ denote the tangent vector space at $x \in X$. Let $c:X \rightarrow Y$ be a differentiable map from a differentiable manifold $X$ into another manifold $Y$. Let ${dc}_x:T_x X \rightarrow T_{c(x)} Y$ denote the differential of $c$ at $x \in X$ and this is a linear map. If the rank of the differential ${dc}_x$ is smaller than the minimum between $\dim X$ and $\dim Y$, then $x$ is a {\it singular point} of $c$.  For a {\it singular} point $x \in X$ of $c$, $c(x)$ is a {\it singular value} of $c$. Let $S(c)$ denote the set of all singular points (the {\it singular set} of $c$). We use "{\it critical}" instead of "singular" for a real-valued function $c:X \rightarrow \mathbb{R}$. Here, we consider smooth maps or maps of the class $C^{\infty}$ as differentiable maps, unless otherwise stated. A {\it diffeomorphism} means a homeomorphism which has no critical point and we can define the notion that two manifolds are {\it diffeomorphic} or equivalently, a notion that a manifold is diffeomorphic to another manifold, under the relation. 
The canonical projection of the Euclidean space ${\mathbb{R}}^k$ into ${\mathbb{R}}^{k_1}$ is denoted
by ${\pi}_{k,k_1}:{\mathbb{R}}^{k} \rightarrow {\mathbb{R}}^{k_1}$ with ${\pi}_{k,k_1}(x)=x_1$ where $x=(x_1,x_2) \in {\mathbb{R}}^{k_1} \times {\mathbb{R}}^{k_2}={\mathbb{R}}^k$ with $k_1, k_2>0$ and $k=k_1+k_2$. The canonical projection of the unit sphere $S^{k-1}$ is its restriction. Hereafter, a {\it real algebraic manifold} means a union of connected components of the zero set of a real polynomial map and a set which is also {\it non-singular}{\rm : }it is defined by the implicit function theorem for the real polynomial map. The space ${\mathbb{R}}^k$, which is also called the {\it $k$-dimensional real affine space}, and the unit sphere $S^{k-1}$ are of simplest examples. {\it Real algebraic} maps are the compositions of the canonical embeddings into the real affine spaces with the canonical projections.\subsection{An SS-region: a region surrounded by circles }
A {\it circle} ({\it centered at} $x_0 \in {\mathbb{R}}^2$) means a real algebraic manifold of the form $\{x \mid ||x-x_0||=r\}$ with its radius $r>0$. 
Hereafter, a {\it region surrounded by circles} or an {\it SS-region} means a pair $(D,\{S_j\}_{j=1}^l)$, where $D$ is a connected open set $D \subset {\mathbb{R}}^2$ with the boundary ${\overline{D}}^{{\mathbb{R}}^2}-D$ represented as a subset of the union ${\bigcup}_{j=1} S_j$  of $l>0$ circles $\{S_j\}_{j=1}^l$ and ${\overline{D}}^{{\mathbb{R}}^2} \bigcap S_j \neq \emptyset$ for each $1 \leq j \leq l>0$. We also assume that for distinct two circles $S_{j_1}$ and $S_{j_2}$, at each point $p_{j_1,j_2} \in S_{j_1} \bigcap S_{j_2} \bigcap ({\overline{D}}^{{\mathbb{R}}^2}-D)$, we have $T_{p_{j_1,j_2}} S_{j_1} \oplus T_{p_{j_1,j_2}} S_{j_2}=T_{p_{j_1,j_2}} {\mathbb{R}}^2$ and that for distinct three circles $S_{j_1}$, $S_{j_2}$, and $S_{j_3}$,  $S_{j_1} \bigcap S_{j_2} \bigcap S_{j_3} \bigcap ({\overline{D}}^{{\mathbb{R}}^2}-D) = \emptyset$. A {\it singular} point $p$ of an SS-region is either of the following.
\begin{itemize}
\item A point $p$ of ${\overline{D}}^{{\mathbb{R}}^2}-D$ contained in exactly two circles of $\{S_j\}_{j=1}^l$.
\item A point $p$ of ${\overline{D}}^{{\mathbb{R}}^2}-D$ contained in exactly one circle of $\{S_j\}_{j=1}^l$ which is a critical point of the restriction of ${\pi}_{2,1}$ to a sufficiently small open neighborhood $U_p \subset {\overline{D}}^{{\mathbb{R}}^2}-D$ of $p$ considered in ${\overline{D}}^{{\mathbb{R}}^2}-D$, diffeomorphic to $\mathbb{R}$.
\end{itemize}
A {\it graph} means a $1$-dimensional finite and connected polyhedron. Its $0$-cell is a {\it vertex} of it and its $1$-cell is an {\it edge} of it. The {\it degree} of a vertex of a graph means the number of edges incident to this.  A {\it tree} is a graph whose 1st Betti number is $0$. 
The {\it Poincar\'e-Reeb graph} of an SS-region can be defined as follows. 
\begin{itemize}
\item The underlying set consists of all components of ${{\pi}_{2,1}}^{-1}(p) \bigcap {\overline{D}}^{{\mathbb{R}}^2}$ and is defined as a quotient space of ${\overline{D}}^{{\mathbb{R}}^2}$.
\item A point in the set is a vertex if it is a component containing some singular points of the SS-region. 
\end{itemize}
On edge $e$, the restriction of the quotient map to the preimage of $e$ gives the structure of a product bundle by so-called Ehresmann's fibration theorem. For each component containing some singular points of the SS-region, we can find finitely many edges and a sequence in each edge converging to it. \cite[Theorem 3.1]{saeki} is a general theorem guaranteeing this where we do not need to understand the theorem. This is motivated by the study \cite{bodinpopescupampusorea}, where regions surrounded by mutually disjoint real algebraic curves in ${\mathbb{R}}^2$ are considered. In addition, as their main result, for a graph embedded in a generic way, such regions are reconstructed.
\subsection{The content of the present paper and our main result.}
In the next section, we introduce two kinds of operations obtaining a new SS-region from a given region by adding a circle and making the region smaller. Our main result is the following. We exhibit again as Main Theorem \ref{mthm:2}. We also give related remarks on reconstruction of {\it real algebraic} maps locally like moment maps.
\begin{MainThm}
\label{mthm:1}
We have a tree $T$ obtained in the following way inductively as the Poincar\'e-Reeb graph of an SS-region obtained by applying finitely many times of the two kinds of operations starting from $(D^2,\{S^1\})$ with $j_0:=0$. 
\begin{itemize}
\item Prepare a graph $G_0$ with exactly one edge. Add $2k_{G_0}$ vertices in the edge where $k_{G_0}$ is an integer.
\item For each vertex $v$ which is of degree $2$ in the graph $G_{j_0}$, we prepare a graph $G_{j_0,v}$ with exactly one edge identify a vertex $v_0$ of the graph with two vertices and $v$, and add $2k_{G_{j_0,v}}$ vertices in the edge of $G_{j_0,v}$ where $k_{G_{j_0,v}}$ is an integer. We have a graph $G_{j_0+1}$ and let $G_{j_0}:=G_{j_0+1}$. We do this one after another as vertices of degree $2$ exist in the graph.
\item By adding vertices to edges of $G_{j_0}$ according to the following rule, we have a tree $T${\rm :} for an edge $e$ of $G_{j_0}$ incident to $0 \leq j_{0,e} \leq 2$ vertices of degree $3$, add $2j_e \geq 2j_{0,e}$ vertices with $j_e$ being an integer.
\end{itemize}
\end{MainThm}


\noindent {\bf Conflict of interest.} \\
The author is also a researcher at Osaka Central
Advanced Mathematical Institute (OCAMI researcher), supported by MEXT Promotion of Distinctive Joint Research Center Program JPMXP0723833165. He is not employed there. This is for our studies
and our study also thanks this. \\
\ \\
{\bf Data availability.} \\
No other data are associated to the paper.
   
\section{On Main Theorem \ref{mthm:1}.}
\subsection{Exposition on two kinds of operations adding circles and changing SS-regions.}

We define the projection ${{\pi}_{2,1}}^{\prime}:{\mathbb{R}}^2 \rightarrow \mathbb{R}$ by ${{\pi}_{2,1}}^{\prime}(x_1,x_2)=x_2$. A point $p$ of ${\overline{D}}^{{\mathbb{R}}^2}-D$ contained in exactly one circle of $\{S_j\}_{j=1}^l$ is a {\it horizontal point} of an SS-region $(D,\{S_j\}_{j=1}^l)$ if it is a critical point of the restriction of ${{\pi}_{2,1}}^{\prime}$ to a sufficiently small open neighborhood $U_p \subset {\overline{D}}^{{\mathbb{R}}^2}-D$ of $p$ considered in ${\overline{D}}^{{\mathbb{R}}^2}-D$, diffeomorphic to $\mathbb{R}$.

By considering ${{\pi}_{2,1}}^{\prime}$ instead of ${\pi}_{2,1}$, we can have a graph similar to the Poincar\'e-Reeb graph of the SS-region. We call this the {\it vertical Poincar\'e-Reeb graph} of the region.

We explain the operations. In addition to \cite{kitazawa3}, we also respect \cite{kitazawa4, kitazawa5}, where we do not assume related knowledge including arguments.

\subsubsection{Adding a small circle centered at a point in the boundary of $D$.}
For an SS-region $(D,\{S_j\}_{j=1}^l)$, add an circle $S_{l+1}$ centered at a point $p \in \overline{D}$ satisfying the following.
\begin{itemize}
\item The point $p$ is contained in exactly one circle $S_j$ from $\{S_j\}_{j=1}^l$.
\item The point is not mapped to any vertex of the Poincar\'e-Reeb graph of the SS-region or that of the vertical Poincar\'e-Reeb graph of the SS-region.
\end{itemize}
If the radius of $S_{l+1}$ is sufficiently small and $S_{l+1}$ bounds a closed disk $D_{l+1} \subset {\mathbb{R}}^2$, then, by putting $D^{\prime}:=D \bigcap ({\mathbb{R}}^2-D_{l+1})$, we have another SS-region $(D^{\prime},\{S_j\}_{j=1}^{l+1})$. 
\begin{Def}
This operation is called an {\it addition of a circle for Morse-Bott circle centered arrangements}, an {\it MBC circle addition}, or an {\it MBCC addition} to $(D,\{S_j\}_{j=1}^l)$.
\end{Def}

Hereafter, we only use the name "SSCC addition" for the operation.

\begin{figure}
	
	\includegraphics[height=75mm, width=100mm]{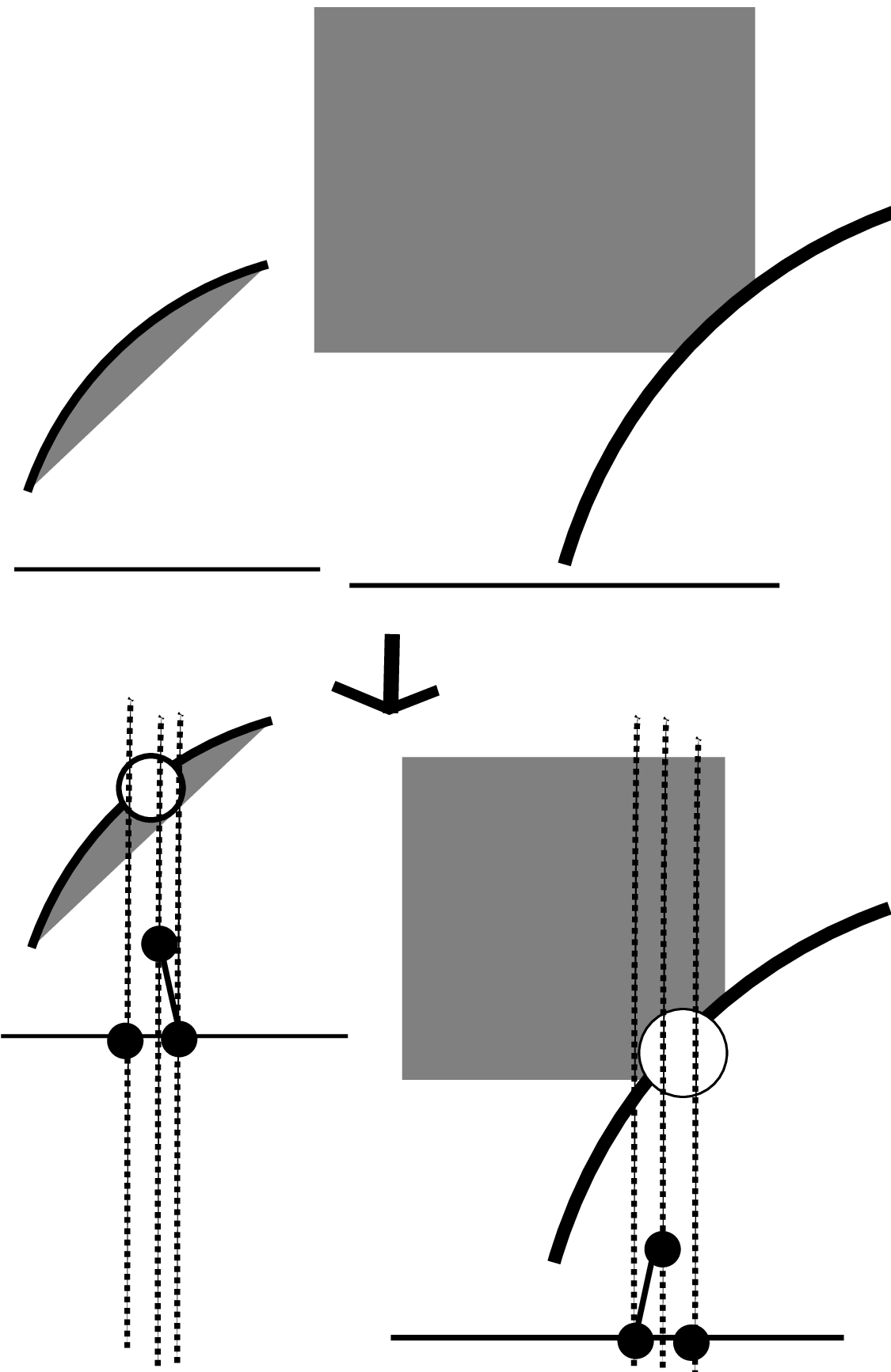}

	\caption{Some cases of MBCC additions are presented. Black colored arcs represent ${\overline{D}}^{{\mathbb{R}}^2}-D$ or ${{\overline{D}}^{\prime}}^{{\mathbb{R}}^2}-D^{\prime}$. Gray colored regions represent part of $D$ or that of $D^{\prime}$.
The Poincar\'e-Reeb graph of $(D,\{S_j\}_{j=1}^{l})$ and that of $(D^{\prime},\{S_j\}_{j=1}^{l+1})$ are also presented for each of the cases. We apply this rule for FIGUREs presented later.}
	\label{fig:1}
\end{figure}

We have the following easily by local observations as in FIGURE \ref{fig:1} with the symmetry. 
\begin{Cor}
\label{cor:1}
If we do an MBCC addition to an SS-region $(D,\{S_j\}_{j=1}^l)$, then the Poincar\'e-Reeb graph of $(D,\{S_j\}_{j=1}^{l})$ changes into that of $(D^{\prime},\{S_j\}_{j=1}^{l+1})$ as follows.
\begin{itemize}
\item An edge $e$ of the Poincar\'e-Reeb graph of $(D,\{S_j\}_{j=1}^{l})$ is chosen and two distinct vertices $v_{e,1}, v_{e,2} \subset e$ are added.
\item A new edge connecting one vertex $v_{e,i}$ of $\{v_{e,1},v_{e,2}\}$ and another new vertex $v_{e,3}$ is added.
\end{itemize}
In this situation, $S_{l+1} \bigcap ({\overline{D}}^{{\mathbb{R}}^2}-D)$ is a three-point discrete set and we have $a_{1,1,S_{{l+1},(D,\{S_j\}_{j=1}^l)}}<a_{1,2,S_{{l+1},(D,\{S_j\}_{j=1}^l)}}<a_{1,3,S_{{l+1},(D,\{S_j\}_{j=1}^l)}}$ as the three distinct values of ${\pi}_{2,1}$ at these three points.
Furthermore, in the present scene, we can choose any edge $e$ of the Poincar\'e-Reeb graph of $(D,\{S_j\}_{j=1}^{l})$ by considering a suitable MBCC addition to $(D,\{S_j\}_{j=1}^{l})$.
\end{Cor}

\subsubsection{Adding a circle passing two distinct and sufficiently close points in the boundary of $D$.}
Hereafter, a {\it straight line} means a non-singular real algebraic set of the form $\{(x_1,x_2) \in {\mathbb{R}}^2 \mid a_1x_1+a_2x_2=0\}$ with $(a_1,a_2) \in {\mathbb{R}}^2-\{(0,0)\}$. A {\it segment} means a subset of a straight line diffeomorphic to $D^1$. A {\it chord of a circle} means a segment connecting two distinct points in the circle, where by the well-known arguments on elementary Euclidean plane geometry, we can have a chord of the circle connecting the points uniquely.

For an SS-region $(D,\{S_j\}_{j=1}^l)$, consider two distinct points $p_1,p_2 \in {\overline{D}}^{{\mathbb{R}}^2}$ satisfying the following.
\begin{itemize}
\item The points $p_1$ and $p_2$ are contained in exactly one circle $S_j$ from $\{S_j\}_{j=1}^l$ and a connected arc in $S_j \bigcap \overline{D}$.
\item Both the points $p_1$ and $p_2$ are mapped to (the interior of) a same edge of the Poincar\'e-Reeb graph of the SS-region and (the interior of) a same edge of the vertical Poincar\'e-Reeb graph of the SS-region.
\end{itemize}
Case A Let the chord $C_{p_1,p_2}$ connecting $p_1$ and $p_2$ be in ${\overline{D}}^{{\mathbb{R}}^2}$. We consider two cases. \\
\ \\
Case A-1\  Suppose that $S_{l+1}$ is a circle passing $p_1$ and $p_2$, bounds the closed disk $D_{l+1}$, and satisfies the following.
\begin{itemize}
\item The circle $S_{l+1}$ is different from any circle of $\{S_j\}_{j=1}^l$.
\item The segment $C_{p_1,p_2}$ is in the closure ${\overline{D_{l,l+1}}}^{{\mathbb{R}}^2}$ of the region $D_{j,l+1}$ surrounded by $S_{j}$ and $S_{l+1}$. 

\end{itemize}
In this case, if $p_1$ and $p_2$ are sufficiently close and the Lebesgue measure of ${\overline{D_{j,l+1}}}^{{\mathbb{R}}^2}$ is sufficiently small, then, as before, by $D^{\prime}:=D \bigcap ({\mathbb{R}}^2-D_{l+1})$, we have another SS-region $(D^{\prime},\{S_j\}_{j=1}^{l+1})$. \\
\ \\
Case A-2\ Suppose that $S_{l+1}$ is a circle passing $p_1$ and $p_2$, and satisfies the following.
\begin{itemize}
\item The circle $S_{l+1}$ is different from any circle of $\{S_j\}_{j=1}^l$.
\item The segment $C_{p_1,p_2}$ is outside the region $D_{j,l+1}$ surrounded by $S_{j}$ and $S_{l+1}$. 

\end{itemize}
In this case, if the Lebesgue measure of ${\overline{D_{j,l+1}}}^{{\mathbb{R}}^2}$ is sufficiently small, then, by $D^{\prime}:=D \bigcap D_{l+1}$, we have another SS-region $(D^{\prime},\{S_j\}_{j=1}^{l+1})$. \\
\ \\
Case B Let the chord $C_{p_1,p_2}$ connecting $p_1$ and $p_2$ be outside ${\overline{D}}^{{\mathbb{R}}^2}$. \\
\ \\ Suppose that $S_{l+1}$ is a circle passing $p_1$ and $p_2$, and $S_{l+1}$ is different from any circle of $\{S_j\}_{j=1}^l$. If the Lebesgue measure of ${\overline{D_{j,l+1}}}^{{\mathbb{R}}^2}$ is sufficiently small, then, as before, by $D^{\prime}:=D \bigcap ({\mathbb{R}}^2-D_{l+1})$, we have another SS-region $(D^{\prime},\{S_j\}_{j=1}^{l+1})$. \\
\begin{Def}
	Each operation presented in the three cases, Case A-1, Case A-2, and Case B, is called an {\it addition of a circle supported by a small chord}, an {\it SSC circle addition}, or an {\it SSCC addition} to $(D,\{S_j\}_{j=1}^l)$. 
\end{Def}
Hereafter, we only use the name "SSCC addition" for the operation.

\begin{figure}
	
	\includegraphics[height=75mm, width=100mm]{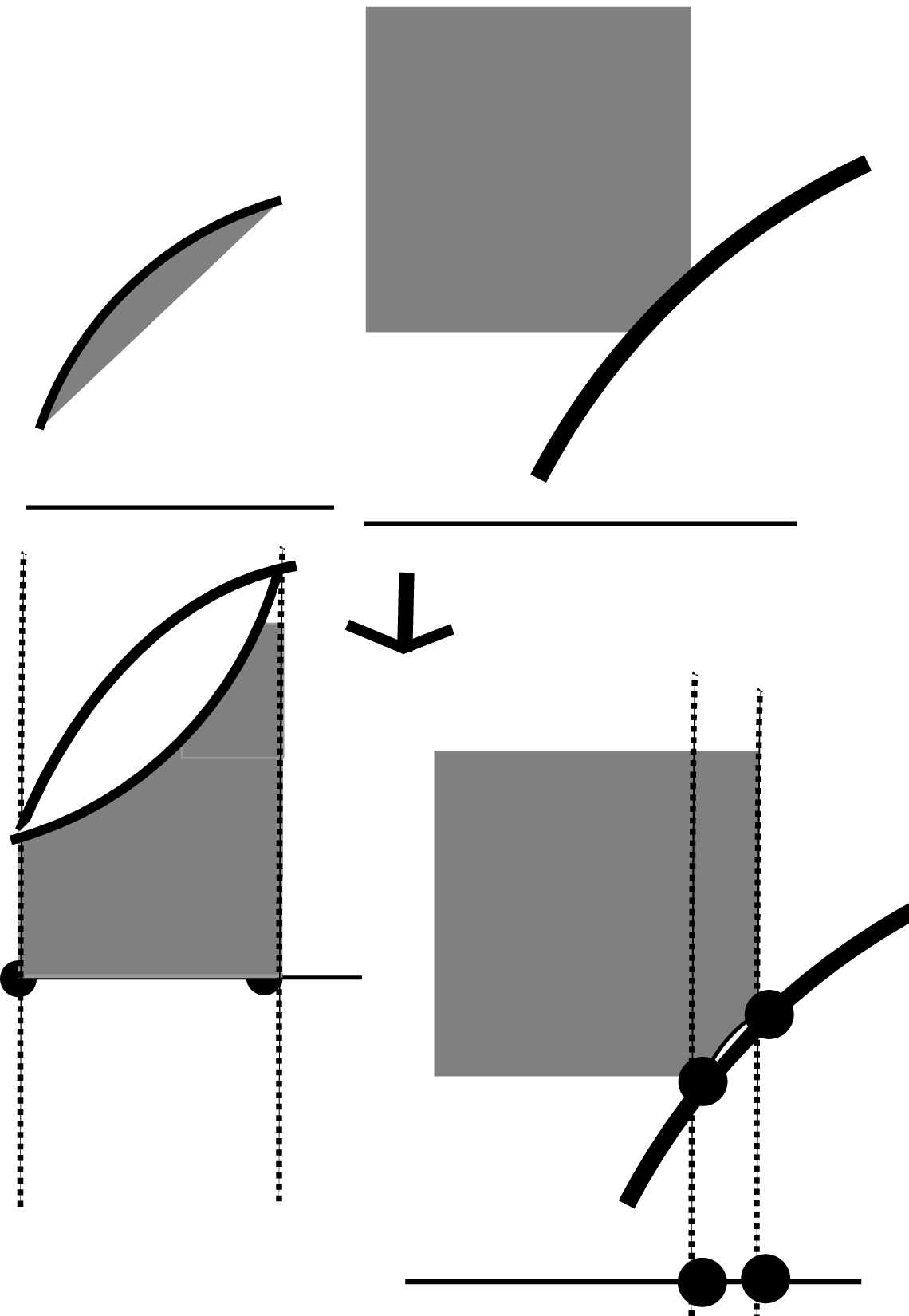}

	\caption{An example for Case A-1 and an example for case B of SSCC additions are presented. The Poincar\'e-Reeb graph of $(D,\{S_j\}_{j=1}^{l})$ and that of $(D^{\prime},\{S_j\}_{j=1}^{l+1})$ are depicted for each case. Case A-2 is omitted.}
	\label{fig:2}
\end{figure}

 We have the following easily by local observations as in FIGURE \ref{fig:2} with the symmetry.
\begin{Cor}
\label{cor:2}
If we do an SSCC addition to an SS-region $(D,\{S_j\}_{j=1}^l)$, then the Poincar\'e-Reeb graph of $(D,\{S_j\}_{j=1}^{l})$ changes into that of $(D^{\prime},\{S_j\}_{j=1}^{l+1})$ as follows.
\begin{itemize}
\item An edge $e$ of the Poincar\'e-Reeb graph of $(D,\{S_j\}_{j=1}^{l})$ is chosen.
\item For the previous edge $e$, two distinct vertices $v_{e,1}, v_{e,2} \subset e$ are added.
\end{itemize}
In this situation, $S_{j+1} \bigcap ({\overline{D}}^{{\mathbb{R}}^2}-D)$ is a two-point discrete set and we have $a_{2,1,S_{{l+1},(D,\{S_j\}_{j=1}^l)}}<a_{2,2,S_{{l+1},(D,\{S_j\}_{j=1}^l)}}$ as the two distinct values of ${\pi}_{2,1}$ at these two points.
Furthermore, in the present scene, we can choose any edge $e$ of the Poincar\'e-Reeb graph of $(D,\{S_j\}_{j=1}^{l})$ by considering a suitable SSCC addition to $(D,\{S_j\}_{j=1}^{l})$.
\end{Cor}

\subsection{Main Theorem \ref{mthm:1} revisited in an improved and revised form.}
We exhibit Main Theorem \ref{mthm:1} again in an improved form.
\begin{MainThm}
\label{mthm:2}
We have a tree $T$ obtained in the following way inductively as the Poincar\'e-Reeb graph of an SS-region obtained by applying finitely many times of the two kinds of operations starting from $(D^2,\{S^1\})$ with $j_0:=0$. 
\begin{itemize}
\item Prepare a graph $G_0$ with exactly one edge. Add $n_{G_0}$ vertices in the edge.
\item For each vertex $v$ which is of degree $2$ in the graph $G_{j_0}$, we prepare a graph $G_{j_0,v}$ with exactly one edge identify a vertex $v_0$ of the graph with two vertices and $v$, and add $n_{G_{j_0,v}}$ vertices in the edge of $G_{j_0,v}$. We have a graph $G_{j_0+1}$ and let $G_{j_0}:=G_{j_0+1}$. We do this one after another as vertices of degree $2$ exist in the graph {\rm (}in other words, we stop the procedure if $n_{G_{j_0,v}}=0$ for $j_0$ and any vertex $v$ of degree $2$ of the graph $G_{j_0}$ {\rm )}.  
\item By adding vertices to edges of $G_{j_0}$ according to the following rule, we have a tree $T${\rm :} for an edge $e$ of $G_{j_0}$ incident to $0 \leq j_{0,e} \leq 2$ vertices of degree $3$, add $n_e$ vertices satisfying the following.
\begin{itemize}
\item In the case $n_{G_0}$ is even and non-zero, for any edge $e$ of $G_{j_0}$ being also a subset of $G_0$, $n_e$ is even and $n_e \geq 2j_{0,e}-1$.  In the case $n_{G_0}$ is odd, for an exactly one arbitrary chosen edge $e_0$ of $G_{j_0}$ being also a subset of $G_0$, $n_{e_0}$ is odd with $n_{e_0} \geq 2j_{0,e_0}$ and for any edge $e$ in the family of all remaining edges of $G_{j_0}$ being subsets of $G_0$, $n_e$ is even and $n_e \geq 2j_{0,e}$. In the case $n_{G_0}=0$, $n_e$ is even and non-negative.
\item In the case $n_{G_{j,v}}$ is even and non-zero, for any edge $e$ of $G_{j_0}$ being also a subset of $G_{j,v}$, relations essentially same as those for $n_{G_0}$ are satisfied.

\end{itemize}
\end{itemize}
\end{MainThm}
\begin{proof}
We can have a tree $G_{j_0}$ with finitely many vertices of degree $2$ added suitably as in the following rule as the Poincar\'e-Reeb graph of an SS-region $(D,\{S_j\}_{j=1}^l)$ by finitely many times of MBCC additions starting from $(D^2,\{S_1:=S^1\})$. Let ${G_{j_0}}^{\prime}$ denote such a graph.
\begin{itemize}
\item There exists a one-to-one correspondence between vertices of degree $2$ of ${G_{j_0}}^{\prime}$ and vertices of degree $3$ of ${G_{j_0}}^{\prime}$.
\item In the correspondence above, the vertex $v_{3,{G_{j_0}}^{\prime}}$ of degree $3$ of the graph ${G_{j_0}}^{\prime}$ and the corresponding vertex $v_{2,{G_{j_0}}^{\prime}}$ of degree $2$ of ${G_{j_0}}^{\prime}$ are adjacent. If $v_{3,{G_{j_0}}^{\prime}}$ is a vertex added in $G_0$, then, $v_{2,{G_{j_0}}^{\prime}}$ is also added in $G_0$. If $v_{3,{G_{j_0}}^{\prime}}$ is a vertex added in $G_{j,v}$, then, $v_{2,{G_{j_0}}^{\prime}}$ is also added in $G_{j,v}$.   
\end{itemize}
This is shown by using Corollary \ref{cor:1} inductively. In addition, we can do in such a way that the following hold.
\begin{itemize}
\item Each circle $S_j$ is centered at a point in $\{x=(x_1,x_2) \in S^1 \mid x_2>0\} \subset S_1=S^1$ and contained in $\{x=(x_1,x_2) \in {\mathbb{R}}^2 \mid x_2>0\}$.
We can do this thanks to the so-called convexity and the shape of the region. Related explicit case is shown in the left part of FIGURE \ref{fig:1}. In short, the newly generated space resulting from an MBCC addition as in this case is regarded as the closure of an new edge of the newly generated graph (the Poincar\'e-Reeb graph of the newly generated SS-region) and the preimage of the natural quotient map of the closure of the region onto the resulting Poincar\'e-Reeb graph contains a smooth connected curve in $\{x=(x_1,x_2) \in S^1 \mid x_2>0\} \subset S_1=S^1$.
\item Either of the following hold for each integer $1 \leq i \leq l-1$, each integer $1 \leq s \leq 2$, each integer $1 \leq t \leq 2$, and
any integer satisfying $1 \leq i^{\prime}< i$. We use $(D_i,\{S_j\}_{j=1}^i)$ for the resulting SS-region after the ($i-1$)-th addition of a circle.
\begin{itemize}
\item $-1<a_{1,s,S_{{i^{\prime}+1},(D_{i^{\prime}},\{S_j\}_{j=1}^{i^{\prime}})}}<a_{1,t,S_{{i+1},(D_i,\{S_j\}_{j=1}^i)}}<a_{1,t+1,S_{{i+1},(D_i,\{S_j\}_{j=1}^i)}}<a_{1,s+1,S_{{i^{\prime}+1},(D_{i^{\prime}},\{S_j\}_{j=1}^{i^{\prime}})}}<1$.
\item $-1<a_{1,s,S_{{i^{\prime}+1},(D_{i^{\prime}},\{S_j\}_{j=1}^{i^{\prime}})}}<a_{1,s+1,S_{{i^{\prime}+1},(D_{i^{\prime}},\{S_j\}_{j=1}^{i^{\prime}})}}<a_{1,t,S_{{i+1},(D_i,\{S_j\}_{j=1}^i)}}<a_{1,t+1,S_{{i+1},(D_i,\{S_j\}_{j=1}^i)}}<1$.
\item $-1<a_{1,t,S_{{i+1},(D,\{S_j\}_{j=1}^i)}}<a_{1,t+1,S_{{i+1},(D,\{S_j\}_{j=1}^i)}}<a_{1,s,S_{{i^{\prime}+1},(D_{i^{\prime}},\{S_j\}_{j=1}^{i^{\prime}})}}<a_{1,s+1,S_{{i^{\prime}+1},(D_{i^{\prime}},\{S_j\}_{j=1}^{i^{\prime}})}}<1$.
\end{itemize}
\item The relation $l=1+n_{G_0}+{\Sigma}_{j=1}^{j_0} ({\Sigma}_{v \in V_{G_{j}}} n_{G_{j-1,v}})$ holds where $V_{G_{j}}$ is the set of all vertices of degree $2$ of the graph $G_{j}$ ($V_{G_{j_0}}=\emptyset$).
\item For distinct two integers $i_{s,1}$ and $i_{s,2}$ satisfying $1 \leq i_{s,1}<i_{s_2} \leq n_{G_0}$, we have $-1<a_{1,1,S_{{i_{s,1}+1},(D_{i_{s,1}},\{S_j\}_{j=1}^{i_{s,1}})}}<a_{1,3,S_{{i_{s,1}+1},(D_{i_{s,1}},\{S_j\}_{j=1}^{i_{s,1}})}}<a_{1,1,S_{{i_{s,2}+1},(D_i,\{S_j\}_{j=1}^{i_{s,2}})}}<a_{1,3,S_{{i_{s,2}+1},(D_i,\{S_j\}_{j=1}^{i_{s,2}})}}<1$.
\item 
 For distinct integers $i_{s,1}$ and $i_{s,2}$ satisfying $1+n_{G_0}+{\Sigma}_{j^{\prime}=1}^{j} ({\Sigma}_{v \in V_{G_{j^{\prime}}}} n_{G_{j^{\prime}-1,v}})<i_{s,1}<i_{s,2} \leq 1+n_{G_0}+{\Sigma}_{j^{\prime}=1}^{j+1} ({\Sigma}_{v \in V_{G_{j^{\prime}}}} n_{G_{j^{\prime}-1,v}})$, we have $-1<a_{1,1,S_{{i_{s,1}+1},(D_{i_{s,1}},\{S_j\}_{j=1}^{i_{s,1}})}}<a_{1,3,S_{{i_{s,1}+1},(D_{i_{s,1}},\{S_j\}_{j=1}^{i_{s,1}})}}<a_{1,1,S_{{i_{s,2}+1},(D_i,\{S_j\}_{j=1}^{i_{s,2}})}}<a_{1,3,S_{{i_{s,2}+1},(D_i,\{S_j\}_{j=1}^{i_{s,2}})}}<1$.
\item Each addition of a circle, the corresponding change of SS regions and their Poincar\'e-Reeb graphs, and each change of the graphs from $G_0$ to $G_{j_0}$, are canonically corresponded.
\end{itemize}

We can also check the following by Corollaries \ref{cor:1} and \ref{cor:2}, our rule, and the situation. We can change each MBCC addition adding $S_{i+1}$ freely by the following operation consisting of an SSCC addition, followed by an MBCC addition, if and only if we need to change. Note that $i_1$ is a suitable integer with $i_1 \geq i$. We also use "${D_{0,j}}$" and "$S_{0,j}$" instead of $D_j$ and $S_j$, respectively, and naturally. FIGURE \ref{fig:3} shows an explicit case.
\begin{itemize}
\item A suitable SSCC addition of Case A-2 or Case B with $a_{2,1,S_{{i_1+1},(D_{0,i_1},\{S_{0,j}\}_{j=1}^{i_1})}}=a_{1,1,S_{{i+1},(D,\{S_j\}_{j=1}^{i})}}$ and $a_{2,2,S_{{i_1+1},(D_{0,i_1},\{S_{0,j}\}_{j=1}^{i_1})}}=a_{1,3,S_{{i+1},(D,\{S_j\}_{j=1}^i)}}$.
\item Another suitable MBCC addition with $a_{2,1,S_{{i_1+1},(D_{0,i_1},\{S_{0,j}\}_{j=1}^{i_1})}}+{\epsilon}_{(D,\{S_j\}_{j=1}^l)}=a_{1,1,S_{{i_1+2},(D_{0,i_1+1},\{S_{0,j}\}_{j=1}^{i_1+1})}}<a_{1,3,S_{{i_1+2},(D_{0,i_1+1},\{S_{0,j}\}_{j=1}^{i_1+1})}}=a_{2,2,S_{{i_1+1},(D_{0,i_1}\{S_{0,j}\}_{j=1}^{i_1})}}-{\epsilon}_{(D,\{S_j\}_{j=1}^l)}$ with ${\epsilon}_{(D,\{S_j\}_{j=1}^l)}>0$ being a suitable and sufficiently small positive number and depending on $(D,\{S_j\}_{j=1}^l)$. The circle $S_{0,i_1+2}$ is centered at a point in $\{x=(x_1,x_2) \in S^1 \mid x_2>0\} \subset S_1=S^1$ and contained in $\{x=(x_1,x_2) \in {\mathbb{R}}^2 \mid x_2>0\}$. 
\end{itemize}  
For the quotient map onto the Poincar\'e-Reeb graph of the given SS-region $(D_{0,i_1},\{S_{0,j}\}_{j=1}^{i_1})$, we consider the preimage of "the edge $e$ of Corollary \ref{cor:2}" chosen in the situation and check the intersection of the preimage and ${\overline{D_{0,i_1}}}^{{\mathbb{R}}^2}$. This consists of two connected smoothly embedded arcs in ${\mathbb{R}}^2$. We choose the place where we put the circle $S_{0,i_1+1}$ as the connected arc different from the connected arc around which we put $S_{i+1}$ originally. After that, we put $S_{0,i_1+2}$ around the original arc. Note that we can regard $S_{0,i_1+2} \subset D_{i+1}$ (remember the rule for the notation $\partial D_{i+1}=S_{i+1}$).
See FIGURE \ref{fig:3} again.
\begin{figure}
	
	\includegraphics[height=75mm, width=100mm]{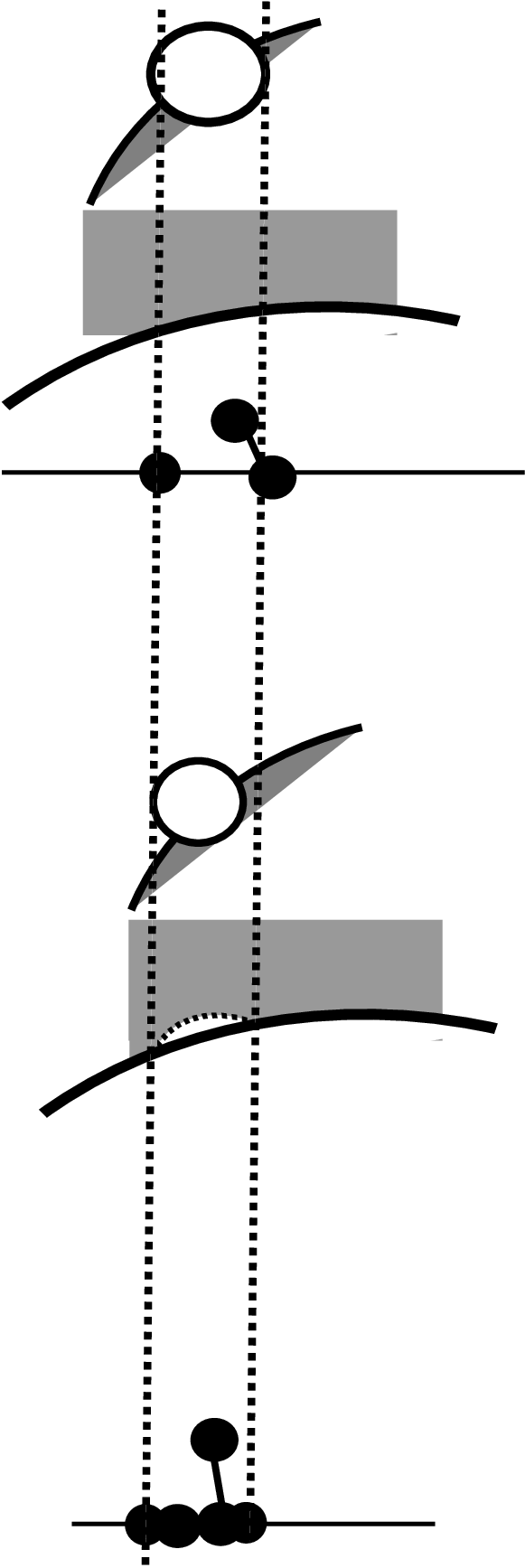}

	\caption{The Poincar\'e-Reeb graph of $(D_{i+1},\{S_{j}\}_{j=1}^{i+1})$ and that of $(D_{0,i_1+2},\{S_{0,j}\}_{j=1}^{i_1+2})$ are depicted for an explicit case of an adjacent pair of an MBCC addition with an SSCC addition. Two additional vertices are added to edges existing already as subsets of the existing graph.}
	\label{fig:3}
\end{figure}

We call this a {\it pair of MBSSCC additions}. This operation increases one vertex to each of the two edges contained already in the existing graph as subsets. This is also presented in FIGURE \ref{fig:3}.  

We compare $G_{j_0}$ and ${G_{j_0}}^{\prime}$. For each edge $e_{G_{j_0}}$ of $G_{j_0}$, $0 \leq j_{0,e_{G_{j_0}}} \leq 2$ is the number of vertices of degree $3$ incident to $e_{G_{j_0}}$.
By doing pairs of MBSSCC additions instead of MBCC additions in some steps, we can have either the following for the resulting graph ${G_{j_0}}^{\prime}$, compared to $G_{j_0}$. We have this fact by the rule for additions of circles and an elementary fact on odd numbers that any odd number is represented as a sum of arbitrary chosen finitely many even numbers and the uniquely defined odd number depending on the finitely many even numbers. 
Let $e_{G_{j_0},{\rm o},G_0}$ denote an arbitrary chosen edge of $G_{j_0}$ contained originally in $G_0$. Let $e_{G_{j_0},{\rm o},G_{j-1,v}}$ denote an arbitrary chosen edge of $G_{j_0}$ contained originally in $G_{j-1.v}$.
\begin{itemize}
\item If $n_{G_0}$ is even and not zero, then the number of additional vertices of ${G_{j_0}}^{\prime}$ contained in each edge $e_{G_{j_0}} \subset G_0$ of $G_{j_0}$ is at most $2j_{0,e_{G_{j_0}}}$ and even.
\item If $n_{G_0}=0$, then the number of additional vertices of ${G_{j_0}}^{\prime}$ contained in the unique edge $e_{G_{j_0}} \subset G_0$ of $G_{j_0}$ is $0$.
\item If $n_{G_0}$ is odd, then the number of additional vertices of ${G_{j_0}}^{\prime}$ contained in each edge $e_{G_{j_0}} \subset G_0$ of $G_{j_0}$ is at most $2j_{0,e_{G_{j_0}}}$ and even except for $e_{G_{j_0},{\rm o},G_0}$, where the number of additional vertices of ${G_{j_0}}^{\prime}$ contained there is odd and at most $2j_{0,e_{G_{j_0},{\rm o},G_0}}-1$.

\end{itemize} 
We can have an essentially same fact for $n_{G_{j-1},v}$ and $e_{G_{j_0},{\rm o},G_{j-1,v}}$.

Last, if we need, we apply SSCC operations for edges of ${G_{j_0}}^{\prime}$ and Corollary \ref{cor:2}. 

This completes the proof.
\end{proof}
\begin{Rem}
For realization of trees as the Poincar\'e-Reeb graphs of SS-regions, see also \cite{kitazawa7} for example.
Note that there the name "SS-region" has not appeared yet. The answers given there are different from the present answers.
\end{Rem}
\subsection{Additional exposition on a moment-like map whose image is the closure of the region $D$ for an SS-region $(D,\{S_j\}_{j=1}^l)$.}
In our setting, by considering a real polynomial $f_j$ of degree $2$ satisfying $S_j=\{x\mid f_j(x)=0\}$ suitably, we can have a situation satisfying $D:={\bigcap}_{j=1}^l \{x \mid f_j(x)>0\}$. We can define a map $l_{(D,\{S_j\}_{j=1}^l)}$ from the finite set $\{S_j\}_{j=1}^l$ onto the set ${\mathbb{N}}_{l^{\prime}}$ of all $l^{\prime} \leq l$ integers from $1$ to $l^{\prime} \leq l$ in such a way that at distinct two $S_{j_1}$ and $S_{j_2}$, the values are distinct.
For each number $i \in {\mathbb{N}}_{l^{\prime}}$, we assign a positive number $d(i)$, where $d$ denotes the map.

We can define the set $M_{(D,\{f_j\}_{j=1}^l)}=\{(x,y)=(x,\{y_i\}_{i=1}^{l^{\prime}} \in {\mathbb{R}}^2 \times {\prod}_{i=1}^{l^{\prime}} {\mathbb{R}}^{d(i)+1}= {\mathbb{R}}^{l^{\prime}+2+{\Sigma}_{i=1}^{l^{\prime}} d(i)} \mid {\prod}_{j \in {l_{(D,\{S_j\}_{j=1}^l)}}^{-1}(i)} (f_j(x))-{||y_i||}^2=0, 1 \leq i \leq l^{\prime}\}$. According to \cite[Main Theorem 1]{kitazawa3}, this is regarded as the zero set of the real polynomial map (we can define naturally) and non-singular. 
By restricting ${\pi}_{l^{\prime}+2+{\Sigma}_{i=1}^{l^{\prime}} d(i),2}$ to $M_{(D,\{f_j\}_{j=1}^l)}$, a real algebraic map is defined. This is locally like a so-called {\it moment map} (in the sense of singularity) and for this, see \cite{buchstaberpanov} and see also \cite{delzant}.

We can also have the function ${\pi}_{l^{\prime}+2+{\Sigma}_{i=1}^{l^{\prime}} d(i),1}$ and this is shown to be a so-called {\it Morse-Bott} function in \cite{kitazawa3}. We do not define a Morse-Bott function rigorously. See \cite{banyagahurtubise} and see also \cite{bott}.

For a smooth function $f:M \rightarrow \mathbb{R}$ on a closed manifold with $S(f)$ being finite, we can define the {\it Reeb graph} $R_f$ whose underlying set is the set of all components of level sets $f^{-1}(t)$ ($t \in \mathbb{R}$) and regarded as the quotient space of $M$ and whose point is a vertex if and only if it is a component containing some critical points of $f$. 
This comes from \cite{saeki} or for more specific cases, see cite{izar} for example. Note that the Reeb graph of a smooth function has been fundamental and strong tools in theory of Morse functions and some functions of certain classes generalizing the class of Morse(-Bott) functions, since the birth of related theory \cite{reeb}.

We can also have the fact that the Reeb digraph $R_f$ and the Poincar\'e-Reeb graph of $(D,\{S_j\}_{j=1}^l)$ are isomorphic.  
 This is related to one of important studies in singularity theory of differentiable maps and applications to algebraic topology and differential topology of manifolds. It is a natural question whether a graph is realized as the Reeb graph of a smooth function (on a closed manifold) of a certain nice class. For this, see \cite{sharko}, a pioneering study on reconstructing nice smooth functions on closed surfaces whose Reeb graphs are as prescribed and whose critical points have certain elementary forms, followed by \cite{masumotosaeki, michalak}. The author has contributed to this by respecting prescribed level sets with no critical point (of the function).

\end{document}